\documentclass[a4paper,leqno]{amsart}
\usepackage{amssymb, amsmath, amscd}
\usepackage{color}

\def\to{\longrightarrow}

\def\cwedge{\bigcirc\kern-1.07em\wedge\ }

\newtheorem{thm}{Theorem}[section]
\newtheorem{lemma}[thm]{Lemma}

\newtheorem{remark}{Remark}

\numberwithin{equation}{section}

\makeatother
\begin{document}
\title[Orthogonal almost complex structures]{Orthogonal almost complex structures\\ on the Riemannian
products\\ of even-dimensional round spheres}
\author{Yunhee Euh and Kouei Sekigawa}
\address{Y. Euh : Department of Mathematics, Sungkyunkwan University, Suwon
440-746, Korea} \email{prettyfish@skku.edu}
\address{K. Sekigawa : Department of Mathematics,
%    Faculty of Science,
    Niigata University,
    Niigata 950-2181, Japan}
\email{sekigawa@math.sc.niigata-u.ac.jp}
\begin{abstract}
We discuss the integrability of orthogonal almost complex structures
on Riemannian products of even-dimensional round spheres and give a
partial answer to the question raised by E. Calabi concerning the
existence of complex structures on a product manifold of a round
2-sphere and a round 4-sphere.\\
\noindent {\it Mathematics Subsect Classification (2010)} : {53C15, 53C21, 53C30}\\
{\it Keywords} : {orthogonal complex structure, even-dimensional
sphere, curvature identity, Ricci $*$-tensor}
\end{abstract}

\maketitle

%%%%%%%%%%%%%%%%%%%%%%%%%%%%%%%%%%%%%%%%%%%%%%%%%%%%%%%%%%%%%%%%%%
\section{Introduction}
%%%%%%%%%%%%%%%%%%%%%%%%%%%%%%%%%%%%%%%%%%%%%%%%%%%%%%%%%%%%%%%%%%

It is well-known that a $2n$-dimensional sphere $S^{2n}$ admits an
almost complex structure if and only if $n=1$ or $3$, and any almost
complex structure on $S^{2}$ is integrable and also the complex
structure on $S^{2}$ is unique with respect to the conformal
structure on it. A 2-dimensional sphere $S^2$ equipped with this
complex structure is biholomorphic to a complex projective line
$\mathbb{C}P_1$. However, contrary to this, it is a long-standing
open problem whether $S^6$ admits an integrable almost complex
structure (namely, complex structure) or not. Lebrun \cite{Le} gave
a partial answer to this problem, that is, proved that any
orthogonal almost complex structure on a round 6-sphere is never
integrable (see also \cite{SV}, Corollary 5.2). On one hand,
Sutherland proved that a connected product of even-dimensional
spheres admits an almost complex structure if and only if it is a
product of copies of $S^2$, $S^6$ and $S^2\times S^4$ under more
general setting (\cite{Su}, Theorem 3.1). In \cite{Ca}, Calabi
raised the problem whether the product manifold $V^2\times S^4$
($V^2$ is any closed, orientable surface) can admit an integrable
almost complex structure or not. In the present note, we discuss the
integrability of orthogonal almost complex structures on a
Riemannian product of round 2-spheres, 6-spheres and Riemannian
product manifolds of a round 2-sphere and a round 4-sphere, and
prove the following.\vskip0.3cm

\noindent{\bf Theorem A.} {\it An orthogonal almost complex
structure on a Riemannian product of round 2-spheres, round
6-spheres, and Riemannian product manifolds of a round 2-sphere and
round a 4-sphere is integrable if and only if it is the product of
{{the}} canonical {{complex}} structures on {{round 2-spheres.}}}
\vskip0.3cm

\begin{remark}\label{remark1}
 % We may give an orthogonal complex structure on any
%Riemannian product of round 2-spheres by the natural way.
Let $M$ be any Riemannian product of round 2-spheres. Then, the
product of the canonical complex structures of the round 2-spheres
is necessarily an orthogonal complex structure on $M$.
\end{remark}

\noindent From Theorem A, we have the following partial answer to
the above mentioned problem by Calabi.\vskip0.3cm

\noindent{\bf Corollary B.} {\it Any orthogonal almost complex
structure on a Riemannain product of a round 2-sphere and a round
4-sphere is never integrable. }\vskip0.3cm

\begin{remark}
An explicit example of an orthogonal almost Hermitian structure on a
Riemannian product of a round 2-sphere and a round 4-sphere was
introduced and its geometric property was discussed in \cite{Ha}.
\end{remark}

We denote by $S^m(\kappa)$ an $m$-dimensional round sphere of
positive constant sectional curvature $\kappa$. Throughout the
present paper, we shall mean by
\textit{a round m-sphere} an oriented $m$-dimensional sphere with constant sectional curvature. \\

{{The authors would like to express their thanks to Professor H.
Hashimoto for drawing their attention
 to the present topic of this paper and also to the referee for his
valuable suggestions.}}

%%%%%%%%%%%%%%%%%%%%%%%%%%%%%%%%%%%%%%%%%%%%%%%%%%%%%%%%%%%%%%%%%%%
\section{Preliminaries}
%%%%%%%%%%%%%%%%%%%%%%%%%%%%%%%%%%%%%%%%%%%%%%%%%%%%%%%%%%%%%%%%%%%%%%%%%%%
Let $M=(M,J,<,>)$ be a $2n$-dimensional almost Hermitian manifold.
We denote by $\nabla$ the Levi-Civita connection and $R$ the
curvature tensor of $M$ defined by
    \begin{equation}\label{def:R}
    R(X,Y)Z=[\nabla_X,\nabla_Y]-\nabla_{[X,Y]}Z
    \end{equation}
for $X$, $Y$, $Z\in \mathfrak{X}(M)$, where $\mathfrak{X}(M)$
denotes the Lie algebra of all smooth vector fields on $M$. We
denote the Ricci $*$-tensor of $M$ by $\rho^*$ which is defined by
    \begin{equation}\label{def:rho*}
    \begin{aligned}
    \rho^*(X,Y)=&\text{tr}\big(Z\longmapsto R(X,JZ)JY\big)\\
               =&\frac{1}{2}\text{tr}\big(Z\longmapsto R(X,JY)JZ\big)
    \end{aligned}
    \end{equation}
for $X$, $Y$, $Z\in \mathfrak{X}(M)$. We here note that the Ricci
$*$-tensor $\rho^*$ satisfies the following equality
    \begin{equation}\label{prop:rho*}
    \rho^*(X,Y)=\rho^*(JY,JX)
    \end{equation}
for $X$, $Y\in \mathfrak{X}(M)$. Thus from \eqref{prop:rho*}, we see
that $\rho^*$ is symmetric if and only if $\rho^*$ is $J$-invariant.
We also denote by $N$ the Nijenhuis tensor of the almost complex
structure $J$ defined by
    \begin{equation}\label{def:N}
    N(X,Y)=[JX,JY]-[X,Y]-J[JX,Y]-J[X,JY]
    \end{equation}
for $X$, $Y\in \mathfrak{X}(M)$. It follows from the celebrated
theorem {of} Newlander and Nirenberg \cite{NN}  that the almost
complex structure $J$ is integrable if and only if $N=0$ {{holds
everywhere on $M$.}} %We may
%also observe that the Nijenhuis tensor $N$ can be expressed in terms
%of the Levi-Civita connection as follows:
%    \begin{equation}\label{eq:N}
%    N(X,Y)=(\nabla_{JX}J)Y-(\nabla_{JY}J)X+J(\nabla_YJ)X-J(\nabla_XJ)Y
%    \end{equation}
%for all $X$, $Y\in \mathfrak{X}(M)$. From \eqref{eq:N}, we may
%easily see that the almost complex structure $J$ is integrable if
%and only if
%    \begin{equation}\label{eq:integrable}
%    (\nabla_{JX}J)JY=(\nabla_XJ)Y
%    \end{equation}
%for all $X$, $Y\in \mathfrak{X}(M)$ {\cite{Gr}}.
An almost Hermitian manifold with an integrable almost complex
structure is called a Hermitian manifold.

Now, we set
    \begin{equation}\label{eq:R}
    R(X,Y,Z,W)=<R(X,Y)Z,W>
    \end{equation}
for $X$, $Y$, $Z$, $W\in \mathfrak{X}(M)$. Gray \cite{Gr} proved the
following result which plays an important role in our forthcoming
arguments of the present paper.
    \begin{thm}\label{thm:Gray}
    The curvature tensor $R$ of a Hermitian manifold $M=(M,J,<,>)$
    satisfies the following identity:
        \begin{equation*}\label{Gray_id}
        \begin{aligned}
        R_{WXYZ}&+R_{JWJXJYJZ}-R_{JWJXYZ}-R_{JWXJYZ}\\
        &-R_{JWXYJZ}-R_{WJXJYZ}-R_{WJXYJZ}-R_{WXJYJZ}=0
        \end{aligned}
        \end{equation*}
    for any $W$, $X$, $Y$, $Z\in\mathfrak{X}(M)$.
    \end{thm}

%%%%%%%%%%%%%%%%%%%%%%%%%%%%%%%%%%%%%%%%%%%%%%%%%%%%%%%%%%%%%%%%%%%
\section{Lammas}
%%%%%%%%%%%%%%%%%%%%%%%%%%%%%%%%%%%%%%%%%%%%%%%%%%%%%%%%%%%%%%%%%%%%%%%%%%%

We shall prepare several lemmas prior to the proof of Theorem A.
First of all, we note that orthogonal almost complex structures on
the Riemannian products of even-dimensional round spheres do not
depend on the order of the factors.
 Now, we consider the Riemannian product $M=S^2(\alpha)\times M'$,
where $M'$ is a Riemannian product of round 2-spheres, round
6-spheres and Riemannian product manifolds of a round 2-sphere and a
round 4-sphere.

    \begin{lemma}\label{lemma2}
    Let $J$ be an orthogonal complex structure on $M$.
    Then, $J$ induces a canonical complex structure on $S^2(\alpha)$ and an orthogonal
    almost complex structure on $\{p_1\} \times M'$ for each point $p_1 \in S^2(\alpha)$.
    \end{lemma}
\textit{Proof of Lemma \ref{lemma2}.} We denote by $\pi_1$ and
$\pi_2$ the canonical projections defined by $\pi_1:M\to
S^2(\alpha)$ and $\pi_2:M\to M'$, respectively. We set
    \begin{equation}\label{eq:3.1}
        x_1= d\pi_1(x),\quad x_2=d\pi_2(x)
    \end{equation}
for any $x\in T_pM$, $p=(p_1,p_2)\in S^2(\alpha)\times M'$. {{ The
tangent space $T_pM$ is identified with the orthogonal direct sum of
$T_{p_1}S^2(\alpha)$ and $T_{p_2}M'$ in the natural way.}}
 Let $x$, $y\in T_{p_1}S^2(\alpha)$ with $x\perp y$, $|x|=|y|=1$. Then, we get
    \begin{equation}\label{eq:R(x,y,x,y)}
    R(x,y,x,y)= -\alpha.
    \end{equation}
Here, since $dim$ $S^2(\alpha)=2$, we may set
    \begin{equation}\label{eq:3.4}
    (Jx)_1=<Jx,y>y,\qquad (Jy)_1=<Jy,x>x.
    \end{equation}
 Now, taking account of
\eqref{eq:3.4}, we get further
    \begin{equation}\label{eq:R(Jx,Jy,Jx,Jy)}
    \begin{aligned}
        &R(Jx,Jy,Jx,Jy)\\
             &=R((Jx)_1+(Jx)_2,(Jy)_1+(Jy)_2,(Jy)_1+(Jy)_2,(Jx)_1+(Jx)_2)\\
             &=R((Jx)_1,(Jy)_1,(Jx)_1,(Jy)_1)+R((Jx)_2,(Jy)_2,(Jx)_2,(Jy)_2)\\
             &=-\alpha \big(|(Jx)_1|^2|(Jy)_1|^2-<(Jx)_1,(Jy)_1>^2\big)\\
             &\qquad+R((Jx)_2,(Jy)_2,(Jx)_2,(Jy)_2)\\
              &=-\alpha |(Jx)_1|^2|(Jy)_1|^2+R_2((Jx)_2,(Jy)_2,(Jx)_2,(Jy)_2),
    \end{aligned}
    \end{equation}
where $R_2$ is the curvature tensor of $M'$.
\begin{equation}\label{eq:R(Jx,Jy,Jx,Jy)}
    \begin{aligned}
   R(Jx,Jy,x,y)
    =&R((Jx)_1,(Jy)_1,x,y)\\
    =&<Jx,y><x,Jy>R(y,x,x,y)\\
    =&\alpha <Jx,y><x,Jy>\\
    =&-\alpha <x,Jy>^2,
    \end{aligned}
    \end{equation}
    \begin{equation}\label{eq:R(Jx,y,Jx,y)}
    \begin{aligned}
    R(Jx,y,Jx,y)
    =&R((Jx)_1,y,(Jx)_1,y)\\
    =&<Jx,y>^2R(y,y,y,y)\\
    =&0,\\
    \end{aligned}
    \end{equation}

    \begin{equation}\label{eq:R(Jx,y,x,Jy)}
    \begin{aligned}
    R(Jx,y,x,Jy)
    =&R((Jx)_1,y,x,(Jy)_1)\\
    =&-<Jx,y>^2R(y,y,x,x)\\
    =&0,
    \end{aligned}
    \end{equation}

    \begin{equation}\label{eq:R(x,Jy,x,Jy)}
    \begin{aligned}
    R(x,Jy,x,Jy)
    =&R(x,(Jy)_1,x,(Jy)_1)\\
    =&<Jy,x>^2R(x,x,x,x)\\
    =&0.
    \end{aligned}
    \end{equation}
Thus, from Theorem \ref{thm:Gray} and
\eqref{eq:R(x,y,x,y)}$\sim$\eqref{eq:R(x,Jy,x,Jy)}, we have
   \begin{equation}\label{eq:3.10}
   \begin{aligned}
    0&=R(x,y,x,y)+R(Jx,Jy,Jx,Jy)-2R(Jx,Jy,x,y)\\
    &\quad-R(Jx,yJx,y)-2R(Jx,y,x,Jy)-R(x,Jy,x,Jy)\\
    &=-\alpha \big\{1-|(Jx)_1|^2|(Jy)_1|^2\big\}+R_2((Jx)_2,(Jy)_2,(Jx)_2,(Jy)_2).\\
   \end{aligned}
   \end{equation}
Since $M'$ is non-negatively curved, we see that
    \begin{equation}\label{eq:R_2}
    R_2((Jx)_2,(Jy)_2,(Jx)_2,(Jy)_2)\leq0
    \end{equation}
for all $x$, $y\in T_{p_1}S^2(\alpha)$. Thus, from \eqref{eq:3.10}
and \eqref{eq:3.1}, we see that
    \begin{equation}\label{eq:3.12}
    |(Jx)_1|=1\quad \text{and} \quad|(Jy)_1|=1
    \end{equation}
and hence, $Jx\in T_{p_1}S^2(\alpha)$ and $Jy\in T_{p_1}S^2(\alpha)$
for any orthogonal pair $\{x,y\}$ in $T_{p_1}S^2(\alpha)$. Since
$d\pi_1$ is a linear map from $T_pM$ onto $T_{p_1}S^2(\alpha)$, from
\eqref{eq:3.12}, we may easily see that $Jx\in T_{p_1}S^2(\alpha)$
for all $x\in T_{p_1}S^2(\alpha)$, and hence
$J(T_{p_1}S^2(\alpha))=T_{p_1}S^2(\alpha)$. Therefore we see also
that $J(T_{p_2}M')=T_{p_2}M'$. \hfill\medskip

 Now, for each $p_1\in S^2(\alpha)$, we denote by $J' = J'(p_{1})$ the induced almost complex structure
on $\{p_1\}\times M'$ in the above Lemma \ref{lemma2}. Then we have
the following.
\begin{lemma}\label{lemma3}
The almost complex structure $J'$ is integrable (and hence defines a
complex structure on $\{p_1\}\times M'$).
\end{lemma}
{\it Proof of Lemma \ref{lemma3}.} %It is evident that the
%submanifold $\{p_1\}\times M'$ of $M$ is totally geodesic.  We
%denote by $\nabla'$ the Levi-Civita connection with respect to the
%induced metric on $\{p_1\}\times M'$. Then, by taking account of
%\eqref{eq:integrable}, we have
%    \begin{equation}\label{eq:3.13}
%    \begin{aligned}
%    &(\nabla'_{J'X'}J')J'Y'\\
%    =&-\nabla'_{J'X'}Y'-J'\nabla'_{J'X'}(J'Y')\\
%    =&-\nabla_{J'X'}Y'-J\nabla_{JX'}(JY')\\
%    =&-\nabla_{J'X'}Y'-J(\nabla_{JX'}J)Y'+\nabla_{J'X'}Y'\\
%    =&(\nabla_{JX'}J)JY'\\
%    =&(\nabla_{X'}J)Y'\\
%    =&\nabla_{X'}(JY')- J\nabla_{X'}Y'\\
%    =&\nabla'_{X'}(J'Y')-J'\nabla'_{X'}Y'\\
%    =&(\nabla'_{X'}J')Y'.
%    \end{aligned}
%    \end{equation}
Let $N'$ be the Nijenhuis tensor of the almost complex
structure $J'$. Then, taking account of Lemma \ref{lemma2}, we have
    \begin{equation}\label{eq:3.14}
    \begin{aligned}
    N'(X',Y')&=[J'X',J'Y']-[X',Y']-J'[J'X',Y']-J'[X',J'Y']\\
             &=[JX',JY']-[X',Y']-J'[JX',Y']-J'[X',JY']\\
             &=[JX',JY']-[X',Y']-J[JX',Y']-J[X',JY']\\
             &=N(X',Y')\\
             &=0
    \end{aligned}
    \end{equation}
for all $X'$, $Y' \in \mathfrak{X}(M')$. Therefore, from
\eqref{eq:3.14}, we see that the induced almost complex structure
$J'$ on {{$\{p_1\}\times M'$ is integrable for each $p_1\in
S^2(\alpha)$}}. \hfill$\square$\medskip

From Lemmas \ref{lemma2} and \ref{lemma3}, if $M'$ involves a round
2-sphere as a factor, {{by a suitable reordering of the factors, we
may assume that $M$ is expressed in the form $M'=S^2(\alpha) \times
M''$}}, where $M''$ expressed by the similar form as $M'$. Then,
applying Lemma \ref{lemma3} to $M'$, it follows that the orthogonal
complex structure $J'$ induces a complex structure on $M''$.
Repeating the similar operations, we may assume that $M$ is
expressed in the form $M= M_1 \times M_2$, where $M_1=
S_{1}^2(\alpha_1) \times \cdots \times S_{s}^2(\alpha_s)$
$(0\leq\alpha_1\leq \cdots \leq \alpha_s)$ and $M_2$ does not
involve a round 2-sphere, and further that the orthogonal almost
complex structure $J$ on $M$ induces a canonical orthogonal complex
structure on $M_1 \times \{p_2\}$ for each point $p_2 \in M_2$ and
an orthogonal almost complex structure on $\{p_1\} \times M_2$ for
each point $p_1 \in M_1$, respectively. Thus, taking account of the
result due to Sutherland (\cite{Su}, Theorem 3.1), we have the
following.
  \begin{lemma}\label{lemma4}
  Let $M$ be a Riemannian product of round 2-spheres, round 6-spheres and
  Riemannian product manifolds of a round 2-sphere and a round 4-sphere, and $J$ be an
  orthogonal complex structure on $M$. Then, $M$ takes of the form $M = M' \times M''$ (after suitable reordering of the factors),
  where $M'$ (resp. $M''$) is a Riemannian product of round 2-spheres (resp. a Riemannian product
  of round 6-spheres), and further, $J$ induces a canonical orthogonal complex structure
  on $M' \times \{p''\}$ for each point $p'' \in M''$ and an orthogonal complex structure on $\{p'\} \times M''$ for each point $p' \in M'$, respectively.
  \end{lemma}
Now, we shall show the following.
   \begin{lemma}\label{lemma5}
   Let $M = (M, < , >)$ be the Riemannian product of round 6-spheres $S_{a}^6(\beta_a) = (S^6, < , >_a)$ $(0<\beta_1\leq\beta_2\cdots\leq\beta_t, a = 1,2,...,t)$,
   and $J$ be an orthogonal almost complex structure on $M$. Then, for each point $(p_1,\cdots,p_{a-1},p_{a+1},\cdots,p_t) \in S_{1}^6(\beta_1) \times \cdots \times S_{a-1}^6(\beta_{a-1})
   \times S_{a+1}^6(\beta_{a+1}) \times \cdots $$\times S_{t}^6(\beta_t)$, $J$ induces an orthogonal almost complex structure on $\{(p_1,\cdots,$ $p_{a-1}, p_{a+1},\cdots,p_t)\}$
   $\times S_{a}^6(\beta_a)$.
   \end{lemma}
\textit{Proof of Lemma \ref{lemma5}.} Let $p=(p_1,p_2,\cdots,p_t)\in
M$ ($p_a\in S_{a}^6(\beta_a)$, $a=1,2,\cdots,t$) be any point of $M$
and $\{e(a)_i\}$ ($i=1,2,\cdots,6$) be any orthonormal basis of
$T_{p_{a}}S_{a}^6(\beta_{a})$. We denote by $R_{(a)}$ the curvature
tensor of $S_{a}^6(\beta_a)$. Then, we have
    \begin{equation}\label{eq:1}
    R(x,y)z=R_{(a)}(x,y)z,
    \end{equation}
and
    \begin{equation}\label{eq:2}
    R_{(a)}(x,y)z=\beta_a\big(<y,z>_a x-<x,z>_a y\big)
    \end{equation}
for $x$, $y$, $z\in T_{p_{a}}S_{a}^6(\beta_a)$. Now, we set
    \begin{equation}\label{eq:3}
    Je(a)_i=\sum_{c=1}^{t}\big(\sum_{j=1}^{6}J(a,c)_{ij}e(c)_j\big)
    \end{equation}
for $1\leq i\leq 6$ and $1\leq a\leq t$. Then since
$<Je(a)_i,e(b)_j> =-<e(a)_i,Je(b)_j>$, from \eqref{eq:3}, we have
    \begin{equation*}
    \begin{aligned}
    <Je(a)_i,e(b)_j>
    &=<\sum_c\sum_kJ(a,c)_{ik}e(c)_k,e(b)_j>\\
    &=\sum_c\sum_kJ(a,c)_{ik}\delta_{cb}\delta_{kj}\\
    &=J(a,b)_{ij}
    \end{aligned}
    \end{equation*}
and
    \begin{equation*}
    \begin{aligned}
    <e(a)_i,Je(b)_j>
    &=<e(a)_i,\sum_c\sum_kJ(b,c)_{jk}e(c)_k>\\
    &=\sum_c\sum_kJ(b,c)_{jk}\delta_{ac}\delta_{ik}\\
    &=J(b,a)_{ji},
    \end{aligned}
    \end{equation*}
and hence, we have
    \begin{equation}\label{eq:4}
    J(a,b)_{ij}=-J(b,a)_{ji}
    \end{equation}
for $1\leq a,b \leq t$ and $1\leq i,j\leq6$. On one hand, since
$J^2=-id$, from \eqref{eq:3}, we have
    \begin{equation*}
    \begin{aligned}
        -e(a)_i&=J(Je(a)_i)\\
               &=J\big(\sum_c\sum_jJ(a,c)_{ij}e(c)_j\big)\\
               &=\sum_c\sum_d\sum_{j,k}J(a,c)_{ij}J(c,d)_{jk}e(d)_k,
    \end{aligned}
    \end{equation*}
and hence,
    \begin{equation}\label{eq:5}
    \sum_c\sum_jJ(a,c)_{ij}J(c,d)_{jk}=-\delta_{ik}\delta_{ad}
    \end{equation}
for $1\leq i,k\leq6$ and $1\leq a,d\leq t$. Here, we shall calculate
the components of the Ricci $*$-tensor $\rho^*$. From \eqref{eq:1},
\eqref{eq:3}, \eqref{eq:4} and \eqref{eq:5}, we have
    \begin{equation}\label{eq:6}
    \begin{aligned}
    \rho^*&(e(a)_i,e(a)_j)\\
    &=-\frac{1}{2}\sum_c\sum_kR(e(a)_i,Je(a)_j,e(c)_k,Je(c)_k)\\
    &=-\frac{1}{2}\sum_kR(e(a)_i,Je(a)_j,e(a)_k,Je(a)_k)\\
    &=-\frac{1}{2}\sum_kR_{(a)}\big(e(a)_i,\sum_lJ(a,a)_{jl}e(a)_l,e(a)_k,\sum_uJ(a,a)_{ku}e(a)_u\big)\\
    &=-\frac{1}{2}\sum_{k,l,u}J(a,a)_{jl}J(a,a)_{ku}
    R_{(a)}(e(a)_i,e(a)_l,e(a)_k,e(a)_u)\\
    &=-\frac{\beta_{a}}{2}\sum_{k,l,u}J(a,a)_{jl}J(a,a)_{ku}
    \{\delta_{lk}\delta_{iu}-\delta_{ik}\delta_{lu}\}\\
    &=-\frac{\beta_{a}}{2}\{-\delta_{ji}-\delta_{ji}\}=\beta_a\delta_{ij},
    \end{aligned}
    \end{equation}
    \begin{equation}\label{eq:7}
    \begin{aligned}
    \rho^*&(e(a)_{i},e(b)_{j})\\
    &=-\frac{1}{2}\sum_{c}\sum_{k}R(e(a)_{i},Je(b)_{j},e(c)_{k},Je(c)_{k})\\
    &=-\frac{1}{2}\sum_{k}R(e(a)_{i},\sum_{l}J(b,a)_{jl}e(a)_{l},e(a)_{k},\sum_{u}J(a,a)_{ku}e(a)_{u})\\
    &=-\frac{1}{2}\sum_{k,l,u}J(b,a)_{jl}J(a,a)_{ku}
    R_{(a)}(e(a)_i,e(a)_l,e(a)_k,e(a)_u)\\
    &=-\frac{\beta_{a}}{2}\sum_{k,l,u}J(b,a)_{jl}J(a,a)_{ku}\{\delta_{lk}\delta_{iu}-\delta_{ik}\delta_{lu}\}\\
    &=-\frac{\beta_{a}}{2}\{J(b,a)_{jk}J(a,a)_{ki}-\sum_lJ(b,a)_{jl}J(a,a)_{il}\}\\
    &=-\frac{\beta_{a}}{2}\{-\delta_{ji}\delta_{ba}-\delta_{ji}\delta_{ba}\}=\beta_{a}\delta_{ij}\delta_{ab},
    \end{aligned}
    \end{equation}
    \begin{equation}\label{eq:8}
    \begin{aligned}
    \rho^*&(e(a)_i,Je(a)_j)\\
    &=\frac{1}{2}\sum_c\sum_kR(e(a)_i,e(a)_j,e(c)_k,Je(c)_k)\\
    &=\frac{1}{2}\sum_kR(e(a)_i,e(a)_j,e(a)_k,Je(a)_k)\\
    &=\frac{1}{2}\sum_{k,l}J(a,a)_{kl}R_{(a)}(e(a)_i,e(a)_j,e(a)_k,e(a)_l)\\
    &=\frac{\beta_{a}}{2}\sum_{k,l}J(a,a)_{kl}\{\delta_{jk}\delta_{il}-\delta_{ik}\delta_{jl}\}\\
%    \end{aligned}
%    \end{equation}
%    \begin{equation*}
%    \begin{aligned}
    &=\frac{\beta_{a}}{2}\{J(a,a)_{ji}-J(a,a)_{ij}\}\\
    &=\beta_aJ(a,a)_{ji},
    \end{aligned}
    \end{equation}

    \begin{equation}\label{eq:9}
    \begin{aligned}
    \rho^*&(e(a)_i,Je(b)_j)\\
    &=\frac{1}{2}\sum_c\sum_kR(e(a)_i,e(b)_j,e(c)_k,Je(c)_k)\\
    &=\frac{1}{2}\sum_{c,d}\sum_{k,l}J(c,d)_{kl}R(e(a)_i,e(b)_j,e(c)_k,e(d)_l)\\
    &=-\beta_{a}\delta_{ab}J(a,b)_{ij},
    \end{aligned}
    \end{equation}

    \begin{equation}\label{eq:10}
    \begin{aligned}
    \rho^*&(Je(a)_i,e(a)_j)\\
    &=-\frac{1}{2}\sum_c\sum_kR(Je(a)_i,Je(a)_j,e(c)_k,Je(c)_k)\\
    &=-\frac{1}{2}\sum_c\sum_{k,l,u,v}J(a,c)_{il}J(a,c)_{ju}J(c,c)_{kv}
    R_{(c)}(e(c)_l,e(c)_u,e(c)_k,e(c)_v)\\
    &=-\frac{1}{2}\sum_c\beta_c\sum_{k,l,u,v}J(a,c)_{il}J(a,c)_{ju}J(c,c)_{kv}\{\delta_{uk}\delta_{lv}-\delta_{lk}\delta_{uv}\}\\
    &=-\frac{1}{2}\sum_c\beta_c\big\{\sum_{k,l}J(a,c)_{il}J(a,c)_{jk}J(c,c)_{kl}\\
    &\qquad\qquad\qquad\quad-\sum_{k,u}J(a,c)_{ik}J(a,c)_{ju}J(c,c)_{ku}\big\}\\
    &=-\frac{1}{2}\sum_c\beta_c\big\{-\sum_lJ(a,c)_{il}\delta_{jl}\delta_{ac}+\sum_uJ(a,c)_{ju}\delta_{iu}\delta_{ac}\big\}\\
    &=\frac{1}{2}\beta_aJ(a,a)_{ij}-\frac{1}{2}\beta_aJ(a,a)_{ji}\\
    &=\beta_aJ(a,a)_{ij},
    \end{aligned}
    \end{equation}

    \begin{equation}\label{eq:11}
    \begin{aligned}
    \rho^*&(Je(b)_i,e(a)_j)\\
    &=-\frac{1}{2}\sum_c\sum_kR(Je(b)_i,Je(a)_j,e(c)_k,Je(c)_k)\\
    &=-\frac{1}{2}\sum_c\sum_{k,l,u,v}J(b,c)_{il}J(a,c)_{ju}J(c,c)_{kv}
    R_{(c)}(e(c)_l,e(c)_u,e(c)_k,e(c)_v)\\
    &=-\frac{1}{2}\sum_c\beta_c\sum_{k,l,u,v}J(b,c)_{il}J(a,c)_{ju}J(c,c)_{kv}
    \{\delta_{uk}\delta_{lv}-\delta_{lk}\delta_{uv}\}\\
    &=-\frac{1}{2}\sum_c\beta_c\big\{\sum_{k,l}J(b,c)_{il}J(a,c)_{jk}J(c,c)_{kl}\\
    &\qquad\qquad\qquad\quad-\sum_{k,u}J(b,c)_{ik}J(a,c)_{ju}J(c,c)_{ku}\big\}\\
    &=-\frac{1}{2}\sum_c\beta_c\big\{-\sum_l\delta_{jl}\delta_{ac}J(b,c)_{il}+\sum_u\delta_{iu}\delta_{bc}J(a,c)_{ju}\big\}\\
    &=\frac{1}{2}\beta_aJ(b,a)_{ij}-\frac{1}{2}\beta_aJ(a,b)_{ji}\\
    &=-\beta_aJ(a,b)_{ij}.
    \end{aligned}
    \end{equation}
Thus, from \eqref{eq:6} and \eqref{eq:7}, we see that $\rho^*$ is
symmetric (and hence, $J$-invariant). Further, from \eqref{eq:9} and
\eqref{eq:11}, taking account of the symmetry of $\rho^*$, we have
$J(a,b)_{ij}=0$ for $a\ne b$, and hence
    \begin{equation}\label{eq:12}
    J(T_{p_{a}}S_{a}^6(\beta_a))=T_{p_{a}}S_{a}^6(\beta_a),\quad
    a=1,2,\cdots,t.
    \end{equation}
Therefore, from \eqref{eq:12}, we see that $J$ induce an almost
complex structure on $\{(p_1,\cdots,p_{a-1},$
$p_{a+1},\cdots,p_t)\}\times S_{a}^6(\beta_a)$ for each
$(p_1,\cdots,p_{a-1},$ $p_{a+1},\cdots,p_t)\in
S_{1}^6(\beta_1)\times \cdots \times S_{a-1}^6(\beta_{a-1})\times$ $
S_{a+1}^6(\beta_{a+1})\times\cdots\times S_{t}^6(\beta_t)$.
 \hfill$\square$\medskip

   \begin{lemma}\label{lemma6}
   Any orthogonal almost complex structure on a Riemannian product of round
   6-spheres is never integrable.
   \end{lemma}
\textit{Proof of Lemma \ref{lemma6}.}
   Let $M = (M,< , >)$ be a Riemannian product of round 6-spheres $S_{a}^6(\beta_a)$ $(a = 1,2,\cdots,t)$, and
   assume that $M$ admits an orthogonal complex structure denoted by $J$. Then, taking account of the results in
   \cite{Le}, it suffices to consider the case where $t \geq 2$. From Lemma \ref{lemma5},
   for each point $(p_1,\cdots,p_{t-1}) \in S_{1}^6(\beta_1) \times \cdots \times S_{t-1}^6(\beta_{t-1})$,
   $J$ induces an orthogonal almost complex structure on
   $\{(p_1,\cdots,p_{t-1})\} \times S_{t}^6(\beta_t)$. Then, we may show that the induced orthogonal almost complex structure
   is integrable by slightly modifying the proof of {{Lemmas \ref{lemma2} and \ref{lemma3}}}.
   But, this is a
   contradiction. \hfill$\square$\medskip

\section{Proof of Theorem A}
{{ In this section, we prove Theorem A based on the arguments in
\textsection 3.}} Let $M = (M, < , >)$ be a Riemannian product of
  round 2-spheres, round 6-spheres, and Riemannian product manifolds of a round 2-sphere and a round 4-sphere,
  assume that $M$ admits an orthogonal complex structure. We denote it by $J$.
 Then, from Lemma \ref{lemma4},
 we see that $M$ is of the form $M = M'\times M''$,
  where $M'$ is of the form $M' = S^2(\alpha_1) \times \cdots \times S^2(\alpha_s)$ and $M''$ is of the form
  $M'' = S^6(\beta_1) \times \cdots \times S^6(\beta_t)$, respectively, and further, $J$ induces an orthogonal complex structure
  on $\{p'\} \times M''$ for each point $p' \in M'$. Therefore, from Lemmas \ref{lemma2} and \ref{lemma3}
  and the uniqueness of the canonical complex structure on a round 2-sphere, $J$ is an orthogonal complex structure on $M$. Therefore, taking account of Lemma \ref{lemma2} we see that $J$ is a product of
  the canonical complex structures on these round 2-spheres. The converse is evident by Remark \ref{remark1}.
   This completes the proof of Theorem A.

\end{document}